\documentclass[11pt, leqno]{amsart}
\usepackage{amsfonts,amsmath,amssymb,amsthm}
\usepackage[margin=1in]{geometry}
\usepackage{url}

\begin{document}

\newtheorem{thm}{Theorem}
\newtheorem{lema}{Lemma}
\newtheorem{prop}{Proposition}
\theoremstyle{definition}
\newtheorem*{remark*}{Remark}

\renewcommand{\Box}{\vrule height 4pt width 4 pt depth 0pt}
\newcommand{\lcm}{\mathop{\mathrm{LCM}}}
\newcommand{\sameauthor}{\leavevmode\hbox to 3em{\hrulefill}}
\newcommand{\emdash}{\hspace{0.2em}\textemdash\hspace{0.2em}}

\makeatletter
	\def\section{\@startsection {section}{1}{\z@}{2em plus 0em minus 0em}{1em plus 0.5em minus 1em}{\centering\large\bf}}
\makeatother
 
\title[Hypergeometric Approach To Euler's Constant]{\bf A Hypergeometric Approach,\\ Via Linear Forms Involving Logarithms,\\to Criteria for Irrationality of Euler's Constant}

\author{Jonathan Sondow\\With an Appendix by Sergey Zlobin}
\date{}


\maketitle
\begin{abstract}
Using an integral of a hypergeometric function, we give necessary and sufficient conditions for irrationality of Euler's constant $\gamma$. The proof is by reduction to known irrationality criteria for $\gamma$ involving a Beukers-type double integral. We show that the hypergeometric and double integrals are equal by evaluating them. To do this, we introduce a construction of linear forms in $1$, $\gamma$, and logarithms from Nesterenko-type series of rational functions. In the Appendix, S.\ Zlobin gives a change-of-variables proof that the series and the double integral are equal.
\end{abstract}

\section{Introduction}
In \cite{Sondow_criteria} we gave criteria for irrationality of Euler's constant $\gamma$, which is defined by the limit
\begin{equation*}
	\gamma := \lim_{N\to\infty} (H_N-\log N),
\end{equation*}
where
\begin{equation*}
	H_N := \sum^N_{k=1} \frac{1}{k}
\end{equation*}
is the $N$th harmonic number. The criteria involve a double integral $I_n$ modeled on Beukers' integrals \cite{Beukers} for $\zeta(2)$ and $\zeta(3)$, and the main step in the proof was to show that
\begin{equation} \label{EQ: d_2n}
	d_{2n} I_n \in \mathbb{Z} + \mathbb{Z}\gamma
		+\mathbb{Z}\log (n+1) + \mathbb{Z}\log (n+2) + \cdots + \mathbb{Z}\log (2n),
\end{equation}
where
\begin{equation*}
	d_n := \lcm(1,2,\dots,n)
\end{equation*}
denotes the least common multiple of the first $n$ natural numbers.

Here we define $I_n$ instead as an integral involving a hypergeometric function; we prove that the same criteria hold with this new $I_n$. The proof is by showing that the old and new definitions of $I_n$ are equivalent. (Alternatively, one could give a self-contained proof along the lines of \cite{Sondow_criteria}; the required inequality $I_n~<~2^{-4n}$ follows easily from Lemma~\ref{LEM: H=S} below.) To show the equivalence, we introduce a series modeled on the one Nesterenko used in \cite{Nesterenko_Zeta3} to give a new proof of Ap\'{e}ry's theorem that $\zeta(3)$ is irrational. (We modify Nesterenko's rational function, and where he differentiates in order to go ``up'' from $\zeta(2)$ to $\zeta(3)$, we integrate to go ``down'' to $\gamma$, which one may think of as ``$\zeta(1)$.'') We prove that both versions of $I_n$ are equal to the sum of our series, by evaluating them. In the Appendix, Sergey Zlobin gives a change-of-variables proof that the double integral and the series are equal, without evaluating them.

The chronology of discovery, different from what one might expect from the above, was as follows. After reading Nesterenko's paper \cite{Nesterenko_Zeta3}, we constructed the series and derived irrationality criteria for $\gamma$ from it. Later, Huylebrouck's survey \cite{Huylebrouck} of multiple integrals in irrationality proofs led us to find the double integral, and using it we rederived the criteria. Zudilin's work \cite{Zudilin} gave us the idea to express the series in hypergeometric form. Here we use Thomae's transformation to simplify the hypergeometric function (compare \cite{Sondow_hyper}). We hope that the variety of expressions for $I_n$ will turn out to be useful in determining the arithmetic nature of $\gamma$.

After seeing \cite{Sondow_hyper}, the Hessami Pilehroods \cite{Pilehrood} extended our non-hypergeometric results to generalized Euler constants. In particular, they used our construction (in Section \ref{SEC: Linear Forms} below) of linear forms involving $\gamma$ and logarithms from series.

For an approach to irrationality criteria for $\gamma$ using Pad\'{e} approximations, see Pr\'{e}vost's preprint \cite{Prevost}.

Recently, Zudilin and the author obtained results \cite{Sondow&Zudilin} analogous to those in \cite{Sondow_criteria}, but using $q$-logarithms instead of ordinary logarithms.

\section{Hypergeometric Irrationality Criteria for Euler's constant}
We state the criteria. First recall that the hypergeometric function $_3F_2$ is defined by the series
\begin{equation*}
	{_3F_2}\biggl(
		\begin{matrix}
			a,\, b,\, c  \\
			d,\, e
		\end{matrix}
	\biggm| z\biggr)
		:= \sum_{k=0}^\infty \frac{(a)_k (b)_k (c)_k}{k! (d)_k (e)_k} z^k\!,
\end{equation*}
where $(a)_0 := 1$ and $(a)_k := a(a+1)\dotsb(a+k-1)$ for $k>0$. The only case we need is when $a,b,c,d,e$ are positive real numbers and $z = 1$. In that case, the series converges if $a+b+c<d+e$. Note that a permutation of either the upper parameters $a,b,c$ or the lower parameters $d,e$ does not change the value of the sum.

Now for $n>0$, let $S_n$ be the positive integer
\begin{equation*}
	S_n := \prod^n_{m=1} \prod^{\min(m-1,n-m)}_{k=0}\prod^{n-k}_{j=k+1}
	(n+m)^{\binom{n}{k}^2{\frac{2d_{2n}}{j}}}
\end{equation*}
and let $I_n$ be the ``hypergeometric integral''
\begin{equation*}
	I_n := \!\int^\infty_{n+1} \frac{(n!)^2\,\Gamma (t)}{(2n+1)\Gamma (2n+1+t)} \, 
		{_3F_2}\biggl(
			\begin{matrix}
				n+1,\, n+1,\, 2n+1  \\
				2n+2,\, 2n+1+t
			\end{matrix}
		\biggm| 1\biggr)\, dt,
\end{equation*}
whose convergence follows immediately from the proof of Lemma~\ref{LEM: H=S}.

\begin{thm}[\textbf{Hypergeometric (ir)rationality criteria for }$\boldsymbol{\gamma}$] \label{THM: HIC}
The following statements are equivalent:
\begin{enumerate}
\setlength{\parskip}{0pt}
\setlength{\itemsep}{0pt}
	\item The fractional part of  $\log S_n$ is equal to $d_{2n}I_n$, for some $n>0$.
	\item The assertion is true for all $n$ sufficiently large.
	\item Euler's constant is a rational number.
\end{enumerate}
\end{thm}

After establishing some preliminary results, we give the proof in Section~\ref{SEC: DblInt for I}.

\section{A Series for $\boldsymbol{I_n}$} \label{SEC: Series for I}
We express the integral $I_n$ as a series.

\begin{lema} \label{LEM: H=S}
If $n>0$, then
\begin{equation}
	I_n = \!\!\sum_{v=n+1}^\infty \int_v^\infty\hspace{-.3em} \left(\frac{n!}{t(t+1)\dotsb(t+n)}\right)^{\! 2} \! dt. \label{EQ: H=S}
\end{equation}
\end{lema}
\begin{proof} Fix $n>0$. We claim that
\begin{equation*}
	I_n = \int^\infty_{n+1} \frac{(n!)^2 \, \Gamma(t)^2}{\Gamma(t+n+1)^2} \,\,
	{_3F_2} \biggl( \begin{matrix}
			t,1,t  \\
			t+n+1, t+n+1
		\end{matrix} \biggm|1\biggr)\, dt.
\end{equation*}
To show this, we apply Thomae's transformation \cite[p.\ 14]{Bailey}, \cite[p.\ 104]{Hardy}, \cite{Krattenthaler}
\begin{equation*}
	\frac{\Gamma(a)}{\Gamma(d) \, \Gamma(e)} \, \,
		{_3F_2} \biggl( \begin{matrix}
			a,\, b,\, c \\
			d,\, e
		\end{matrix} \biggm|1\biggr)
	= \frac{\Gamma(s)}{\Gamma(s+b) \, \Gamma(s+c)} \, \,
		{_3F_2} \biggl( \begin{matrix}
			s,\, d-a,\, e-a \\
			s+b,\, s+c
		\end{matrix} \biggm|1\biggr),
\end{equation*}
where $s := d+e-a-b-c$. After permuting the upper parameters in the resulting function ${_3F_2}$, we obtain the integrand in the definition of $I_n$, proving the claim.

Now since
\begin{equation*}
	{_3F_2} \biggl( \begin{matrix}
		t,\, 1,\, t \\
		t+n+1,\, t+n+1
	\end{matrix} \biggm|1\biggr)
	= 1 + \frac{t^2}{(t+n+1)^2} + \frac{t^2(t+1)^2}{(t+n+1)^2(t+n+2)^2} + \dotsc
\end{equation*}
we can use the identity $\Gamma(x+1) = x\Gamma(x)$ to write
\begin{equation*}
	I_n = \int_{n+1}^\infty \sum_{\nu=0}^\infty
			\frac{(n!)^2 \, \Gamma(t+\nu)^2}{\Gamma(t + \nu + n + 1)^2} \, dt
		= \int_{n+1}^\infty \sum_{\nu=0}^\infty R_n(t+\nu) \, dt,
\end{equation*}
where
\begin{equation*}
	R_n(t) := \biggl(\frac{n!}{t(t+1)\cdots(t+n)}\biggr)^{\!\! 2}
\end{equation*}
is the rational function in \eqref{EQ: H=S}. Interchanging integral and summation, and replacing $t$ with $t-\nu$, and $\nu$ by $\nu-n-1$, we arrive at \eqref{EQ: H=S}.
\end{proof}

\section{Constructing Linear Forms in $\boldsymbol{1}$, $\boldsymbol{\gamma}$ and Logarithms From Series} \label{SEC: Linear Forms}
We give a method for constructing linear forms involving $\gamma$ and logarithms from certain series of rational functions.

\begin{prop} \label{PROP: R&BLA}
Fix $n > 0$ and let $R(t)$ be a rational function over $\mathbb{C}$ of the form
\begin{equation}
	R(t) = \sum^n_{k=0}\biggl( \frac{B_{k2}}{(t+k)^2} +\frac {B_{k1}}{t+k}\biggr).  \label{EQ: R(t)}
\end{equation}
If $R(t)= O(t^{-3})$ as $t\rightarrow\infty$, then
\begin{equation*}
	\sum^\infty_{\nu=n+1} \int^\infty_\nu \!\! R(t)\, dt = B\gamma + L - A,
\end{equation*}
where
\begin{equation}
	B := \sum^n_{k=0} B_{k2}, \quad
	L := \sum^n_{m=1}\sum^n_{k=m}B_{k1} \log(n+m), \quad
	A := \sum^n_{k=0} B_{k2} H_{n+k}.
\end{equation}
\end{prop}
\begin{proof}
From \eqref{EQ: R(t)},  for $|t|$ large we obtain an expansion $R(t)=\sum^\infty_{i=1} b_i t^{-i}$\!, with $b_1 = \sum^n_{k=0} B_{k1}$ and $b_2 = \sum^n_{k=0} (B_{k2} - kB_{k1})$. The asymptotic hypothesis implies that $b_1 = b_2 = 0$, so we have the relations
\begin{equation}
	\sum^n_{k=0} B_{k1} =0, \quad
	\sum^n_{k=0} (B_{k2} - kB_{k1})=0.  \tag{5.1, 5.2} \label{EQ: 2part}
\end{equation}
\addtocounter{equation}{1}

In view of (5.1), the sums $\sum^n_{k=0} B_{k1}\log (t+k)$ and $\sum^n_{k=0} B_{k1}\log (1+kt^{-1})$ are equal.  Hence for $N>n$ we have
\begin{equation}
	\sum^N_{\nu=n+1} \int^\infty_\nu \hspace{-0.3em} R(t) \, dt
		= \hspace{-0.5em} \sum^N_{\nu=n+1} \sum^n_{k=0}
			\left( \frac{B_{k2}}{\nu+k} - B_{k1}\log (\nu+k) \right).  \label{EQ: BLA}
\end{equation}
Define $B$, $L$, $A$ by \eqref{EQ: BLA}, and rewrite $L$ as $L= \sum^n_{k=1} \sum^{n+k}_{m=n+1} B_{k1}\log m$. Evidently the double sum in \eqref{EQ: BLA} differs from the expression
\begin{equation}
	BH_N + L-A + \sum^n_{k=1} \, \sum^{N+k}_{m=N+1}
		\hspace{-0.3em}\left( \frac{B_{k2}}{m}-B_{k1}\log m \right)  \label{EQ: BH_N}
\end{equation}
by the quantity $\sum^n_{k=0}\sum^N_{m=n+1} B_{k1}\log m$, which vanishes by (5.1).  Since $n$ is fixed and $k\leq n$, the double sum in \eqref{EQ: BH_N} equals $-\sum^n_{k=0} kB_{k1}\log N + O(N^{-1})$ as $N\rightarrow\infty$. Using (5.2), it follows that the left-hand side of \eqref{EQ: BLA} is equal to $B(H_N-\log N) + L-A+O(N^{-1})$, and we obtain the required
formula by letting $N$ tend to infinity.
\end{proof}

\section{Summing the Series for $\boldsymbol{I_n}$}
Applying Proposition~\ref{PROP: R&BLA} to the rational function $R(t)$, we sum series \eqref{EQ: H=S} for $I_n$.

\begin{lema} \label{LEM: S=LF}
If $n > 0$, then
\begin{equation*}
	I_n = \binom{2n}{n}\gamma + L_n - A_n,
\end{equation*}
where $L_n$ is the linear form in logarithms
\begin{equation*}
	L_n := \sum^n_{m=1} \sum^{\min (m-1, n-m)}_{k=0} \sum^{n-k}_{j=k+1}
		\! \binom{n}{k}^{\hspace{-0.25em} 2} \, \frac {2}{j} \log (n+m)
\end{equation*}
and $A_n$ is the rational number
\begin{equation*}
	A_n := \sum^n_{k=0} \binom{n}{k}^{\hspace{-0.25em} 2} H_{n+k}.
\end{equation*}
Moreover, inclusion \eqref{EQ: d_2n} and $d_{2n} L_n = \log S_n$.
\end{lema}
\begin{proof}
The partial fraction decomposition of the integrand in \eqref{EQ: H=S} is given by the right-hand side of \eqref{EQ: R(t)} with
\begin{align*}
	B_{k2} = (t+k)^2 R_n(t) |_{t=-k} = \binom{n}{k}^{\hspace{-0.25em} 2}
\end{align*}
and
\begin{align*}
	B_{k1} = \frac{d}{dt} \Bigl( (t+k)^2 R_n(t)\Bigr) \hspace{-0.3em}\Bigm|_{t=-k}
		= 2\binom{n}{k}^{\hspace{-0.25em} 2} \! \left(H_k - H_{n-k}\right),
\end{align*}
where $H_0 = 0$. Using the relations $\sum^n_{k=0} \binom{n}{k}^{\! 2}\! = \binom{2n}{n}$ and
$B_{n-k, 1} =-B_{k1}$, the result follows from Proposition~\ref{PROP: R&BLA} and the definition of $S_n$.
\end{proof}

\section{A Double Integral For $\boldsymbol{I_n}$ and Proof of the Criteria} \label{SEC: DblInt for I}
We obtain another representation of $I_n$, as a double integral, and prove Theorem \ref{THM: HIC}.

\begin{lema} \label{LEM: S=IntInt}
For $n>0$, the following equality holds:
\begin{equation}
	\sum^\infty_{\nu=n+1} \int^\infty_\nu \biggl(\frac{n!}{t(t+1)\cdots(t+n)}\biggr)^{\!2} dt
		= \iint\limits_{[0,1]^2} \frac{\bigl( x(1-x) y(1-y)\bigr)^n}{(1-xy)(-\log xy)}\, dx \, dy.
		\label{EQ: S=IntInt}
\end{equation}
\end{lema}
\begin{proof} By Lemmas \ref{LEM: H=S}, \ref{LEM: S=LF}, and \cite[(8)]{Sondow_criteria}, the series and the double integral, respectively, are both equal to $\binom{2n}{n}\gamma + L_n - A_n$.
\end{proof}

\begin{proof}[Proof of Theorem \ref{THM: HIC}]
This follows immediately from Lemmas \ref{LEM: H=S} and \ref{LEM: S=IntInt} together with the main result of \cite{Sondow_criteria}, which is the same as Theorem~\ref{THM: HIC}, except that in \cite{Sondow_criteria} we defined $I_n$ to be the double integral in \eqref{EQ: S=IntInt}.
\end{proof}

\begin{remark*}
There exist representations of many constants as double integrals of the same shape as the one in \eqref{EQ: S=IntInt}\emdash see \cite{G&S}, \cite{Sondow_dblint}, \cite{Sondow_faster}, \cite{Sondow&Hadjicostas}.
\end{remark*}

\section*{Acknowledgment}
I am grateful to Wadim Zudilin for several suggestions.

\section*{Appendix by Sergey Zlobin}
In this appendix we prove the statement of Lemma~\ref{LEM: H=S} without expanding integrals to linear forms. First, we develop $1/(1-xy)$ in a geometric series
\begin{align*}
	J_n &:= \iint\limits_{[0,1]^2} \frac{x^{n}(1-x)^{n}y^{n}(1-y)^{n}}{(1-xy) (-\log xy)} \,dx \, dy \\[0.3em]
		&= \sum_{k=0}^{\infty} \hspace{0.4em} \iint\limits_{[0,1]^2}
			\frac{x^{n+k} (1-x)^{n}y^{n+k} (1-y)^{n}}{-\log xy}\, dx \, dy.
\end{align*}
(To justify the interchange of the sum and the double integral, one can expand $1/(1-xy)$ in a finite sum with remainder and make the same estimations as in \cite{Sondow_criteria}.) Further, we substitute
\begin{equation*}
	-\frac{(xy)^k}{\log xy} = \int_{k}^{\infty}\!\! (xy)^t dt
\end{equation*}
and obtain
\begin{align*}
	J_n
	&= \sum_{k=0}^{\infty} \hspace{0.4em}\iint\limits_{[0,1]^2}
		\biggl( \int_{k}^{\infty} \!\!\! x^{n+t} (1-x)^{n} y^{n+t} (1-y)^{n}\, dt \biggr)\, dx \, dy \\
	&= \sum_{k=0}^{\infty} \int_{k}^{\infty} \biggl(\,\iint\limits_{[0,1]^2}x^{n+t} (1-x)^{n} y^{n+t}(1-y)^{n}\,
dx \, dy\biggr) \,dt,
\end{align*}
where we can change the order of integration because the integrand is nonnegative and all the integrals converge. Since
\begin{equation*}
	\int_0^1\! u^{n+t} (1-u)^{n} \,du=\frac{n!}{(t+n+1)(t+n+2)\cdots(t+2n+1)}
\end{equation*}
we have
\begin{align*}
	J_n
	&= \sum_{k=0}^{\infty} \int_{k}^{\infty} \!
		\biggl(\frac{n!}{(t+n+1)(t+n+2)\cdots(t+2n+1)^2}\biggr)^{\! 2} dt \\[0.2em]
	&= \!\!\sum_{k=n+1}^{\infty} \int_{k}^{\infty} \! \biggl(\frac{n!}{t(t+1)\cdots(t+n)}\biggr)^{\! 2} dt,
\end{align*}
and we get the desired identity. 

The same method can be applied to prove that (minus) the series Nesterenko used in \cite{Nesterenko_Zeta3} is equal to Beukers' triple integral in \cite{Beukers}. Another proof of that fact is given in \cite{Nesterenko_Integrals} and uses an identity with a complex integral.

\newpage
\setlength{\parindent}{0em}

\vspace{4em}
\hspace{2em}
\small
\begin{minipage}{3in}
\textit{J.\ Sondow\\
209 West 97th Street\\
New York, NY~~10025 USA\\
\url{jsondow@alumni.princeton.edu}}
\end{minipage}
\begin{minipage}{3in}
\textit{S.\ Zlobin\\
Faculty of Mechanics and Mathematics\\
Moscow State University\\
Leninskie Gory, Moscow 119899 RUSSIA\\
\url{sirg_zlobin@mail.ru}}
\end{minipage}

\end{document}